\documentclass{amsart}
\usepackage[dvips]{geometry}

\usepackage{amsfonts}
\usepackage{verbatim}
\usepackage{stmaryrd}
\usepackage{amsmath}
\usepackage[mathscr]{eucal}
\usepackage{color}

\newenvironment{pf}{\noindent{\it Proof.} }{\fbox{}\\}
\newtheorem{thm}{Theorem}
\newtheorem{cor}{Corollary}
\newtheorem{prop}{Proposition}
\newtheorem{lem}{Lemma}
\newtheorem{defi}{Definition}
\newtheorem{rem}{Remark}
\newcommand{\1}{\mathbf{1}}
\newcommand{\+}{\dag}
\newcommand{\club}{$\clubsuit$}
\renewcommand{\d}{\mathrm{d}}
\newcommand{\e}{\mathrm{e}}
\renewcommand{\i}{\mathrm{i}}

\renewcommand{\Im}{\mathrm{Im}}
\renewcommand{\r}{\mathfrak{r}}
\newcommand{\<}{\langle}
\renewcommand{\>}{\rangle}
\newcommand{\dc}{\overline{\partial}}
\title{Results on the Spectral Schwartz Distribution}
\author{Wilhelm von Waldenfels\\
Universit\"at Heidelberg}
\date{}

\begin{document} 

\begin{abstract}
 The resolvent of an operator in a Banach space  is defined on an open subset of the complex plane and is holomorphic. 
 It obeys the resolvent equation. A generalization of this equation to Schwartz distributions is defined and a Schwartz distribution, which
 satisfies that equation is called a resolvent distribution. Its restriction to the subset, where it is continuous, is the usual resolvent function.
 Its complex conjugate derivative is ,but a factor, the spectral Schwartz distribution ,
which is carried by a subset of the spectral set of the operator. The spectral distribution yields a spectral decomposition.
The spectral distribution of a matrix and a unitary operator are given. If the the 
operator is a self-adjoint operator on a Hilbert space, the spectral distribution is the derivative of the spectral family. We calculate the
spectral distribution of the multiplication operator and some rank one perturbations. These operators are not necessarily self adjoint
and may have discrete real or imaginary eigenvalues or a nontrivial Jordan decomposition.
\end{abstract}

\maketitle
infinite dimensional spaces.
\section{Definition and Basic Properties}
\subsection{Introduction}In a study of radiation transfer Garyi V. Efimov, Rainer Wehrse  and the author 
\cite{EvW} developed a method, how to calculate the
spectral decomposition of a  multiplication operator in $\mathbb{R}^n$ perturbed by a  rank one operator
. By Krein's formula it is possible to calculate the resolvent
$R(z)$ and then one investigates the singularities of the  function $z \mapsto\langle f_1| R(z)|f_2 \rangle$, where $f_1,f_2$ are 
test-functions. In \cite[chapter 3]{vW} the complex derivative derivative $ 1/\pi\dc \langle f_1| R(z)|f_2 \rangle$ is called the bracket
$\langle f_1|M(z)|f_2\rangle$ of the 
{\em spectral Schwartz distribution } $M(z)$, which  therefor was only scalarly defined.
 There was developed  a rudimentary theory and calculated four explicit examples \cite[chapter 4]
{vW}, originated from the theory of 	quantum stochastic processes.
 Inspecting the examples it turned out, that the resolvent it self, not only
the brackets between test functions can be extended to an operator valued
Schwartz distribution and under this much stronger regularity condition
  a  deeper theory of the spectral Schwartz  distribution
can be established. 
 We calculate the spectrum of the multiplication operator,  and two of  its perturbations by a rank one operator. The examples
are not necessarily self adjoint. In the second example we have  imaginary  eigenvalues and a nontrivial Jordan decomposition. 
\subsection{Schwartz distributions}
We recall the definition of a distribution given byf L.Schwartz   \cite{Schwartz1}. A distribution on on an open set
 $G \subset\mathbb{R}^n$ is a linear functional $T$
 on the space of {\em test functions} $\mathcal{D}(G)$,
that is the space
of infinitely differentiable functions of compact support (  or $C_c^\infty$-functions) with support in $G$.   The functional $T$ 
 is continuous in the following way: if the functions $\varphi_n\in \mathcal{D}(G)$ have their support
in a common compact set and if they  and all their derivatives converge to 0 uniformly,
 then the $T(\varphi_n)\to 0$. The space of distributions on
$\mathbb{R}^n$ is denoted by $\mathcal{D}'(\mathbb{R}^n)$.
We write often
$$ T(\varphi)=\int T(x) \varphi(x) \d x.$$
The integral notation is the usual notation of physicists and  has been adopted by Schwartz with minor modifications in his articles
 on vector-valued distributions \cite{Schwartz2}\cite{Schwartz3}. It emphasizes the fact, that a distribution can be considered
as a {\em generalized function}. This is the notation in the Russian literature \cite{Gelfand1}.  Similar is
 de Rham's formulation \cite{Rham} : $T$ can be considered as a current of degree 0 applied to a form
of degree $n$ namely
$ \varphi(x_1,\cdots,x_n) \d x_1 \cdots \d x _n$. We use variable transforms of distributions accordingly.

In order to include  Dirac's original ideas, L.Schwartz introduced the {\em distribution kernels} \cite{Schwartz2}.
A distribution kernel $K(x,y)$ on $\mathbb{R}^m \times \mathbb{R}^n$ is a distribution on that space.
 If $T$ is a distribution on $\mathbb{R}^n$ one may associate to it a kernel $T(x-y)$ on
$\mathbb{R}^n \times \mathbb{R}^n$  by
$$ \iint T(x-y) \varphi(x) \psi(y) \d x \d y = \int \d x T(x)\varphi(x)\int \d y \psi(y-x) = \int \d x(T*\psi)(x)\varphi(x) ,$$
where $T*\psi$  denotes  the convolution.  We have often to deal with the so called
$\delta$-function given by $\delta(\varphi)=\varphi(0)$ and the generalized function
$ \mathcal{P}/x = (\d/\d x)\ln|x|$, where $\mathcal{P}$ stands for principal value.We have
$$ \iint \delta(x-y) \varphi(x,y) \d x \d y = \int \d x \varphi(x,x).$$So
$ \delta(x-y)=\delta(y-x)$.

If $S(x,y)$ and $T(y,z)$ are two kernels one may form a distribution of three variables 
$S(x,y)T(y,z)$, if this exists. If the kernels are $S(x-y)$ and $T(y-z)$ then $S(x-y)T(y-z)$ exists and  is given by
$$ \iiint \d x \d y \d z S(x-y)T(y-z) \varphi(x,y,z)= \iiint \d u \d v \d y S(u)T(v) \varphi(u+y,y,y-v).$$
We obtain for any distribution $T$ the relation
\begin{equation}\label{eq1}\delta(x-y)T(y-z)= \delta(x-y)T(x-z). \end{equation}
We have the formulae ( for the second equation see \cite[p.74]{vW})
\begin{align}\label{eq2} &\iiint \delta(x-y)\delta(y-z) \varphi(x,y,z)\d x \d y \d z=  \int \varphi(x,x,x) \d x \\&\label{eq3}
\frac{\mathcal{P}}{x-\omega}\frac{\mathcal{P}}{y-\omega}=
\frac{1}{y-x}\left(\frac{\mathcal{P}}{x-\omega} -\frac{\mathcal{P}}{y-\omega}
\right) + \pi^2 \delta(x-\omega)\delta(y-\omega) \end{align}

\subsection{Notion of the Spectral Distribution}
   Assume we have a Banach space $V$ and
     denote by $L(V)$ the space of all
bounded linear operators from $V$ to $V$ provided with the usual operator norm, Assume
a subspace $D \subset V$ and an operator $A:D\to V$. Then the complex plane splits into 
two sets, the {\em resolvent set} , where the  {\em resolvent} $R(z)= (z-A)^{-1}$
exists, and the {\em spectral set}, where it does not. The resolvent set is open, the spectral set 
is closed. In the resolvent set the resolvent obeys one of the two equivalent  {\em resolvent equations} 
 \begin{align}\label{eq4}&R(z_1)-R(z_2)=(z_2-z_1)R(z_1)R(z_2)\\&\label{eq5}
 R(z_1)R(z_2)= 1/(z_2-z_1)(R(z_1)-R(z_2).\end{align}

Go the other way round and assume
   an open set $G\subset\mathbb{C}$ and
 a function $R(z):G\to L(V)$ satisfying the resolvent equation.
The function  $R(z)$
      is \emph {holomorphic} in $G$.
The subspace $ D= R(z)V$ is a subset independent of $z\in G$. If
$R(z_0)$ is injective for one $z_0\in G$, then $R(z)$ is injective for all $z\in
G$ and there exists a mapping $A:D\to V$ such that
 \begin{align*} (z-A)R(z)f&=f
\text{ for } f\in V & R(z)(z-A)f&=f \text{ for } f\in D \end{align*}
  The operator $A$ is closed and $R(z)$ is the resolvent of $A$.\cite{Hille}
 There exists not always such an operator corresponding to a resolvent , e.g. the function $R(z)=0$ for all $z$ fulfills the 
 resolvent equation and $R(z)$ is surely not injective.

We use the second equation  (\ref{eq5}) and formulate:
\begin{defi} Assume a Schwartz distribution $\varphi \mapsto	R(\varphi)$ on an  open set $G\subset \mathbb{C}$
, which satisfies the equation
\begin{multline*}R(\varphi_1)R(\varphi_2)= \int \d^2 z_1R(z_1) \varphi_1(z_1)
\int \d^2 z_2 \varphi_2 (z_2) /(z_2-z_1)\\
- \int \d^2 z_2R(z_2) \varphi_2(z_2)
\int \d^2 z_1 \varphi_1 (z_1) /(z_2-z_1)\end{multline*}
then we say, that $R(\varphi)$ satisfies the {\em distribution resolvent equation} and call $R(\varphi)$ a {\em  resolvent distribution}.
Here $\d^2 z = \d x \d y$ is the surface elament on the complex plane for $z=x + \i y$.
\end{defi}

 If $R(z)$ obeys the resolvent equation in $G$
 it might be, that there is Schwartz distribution on  an open set $G' \supset G$
extending $R(z)$ and obeying
 the distribution resolvent equation.
Assume a distribution $R(z)$ satisfying the distribution resolvent equation , whose restriction
on an open subset of the complex plane equals a continuous function, then the restriction
satisfies the usual resolvent equation and is holomorphic. We call the set, where the distribution $R(z)$ is  not continuous, the
{\em spectral set}. It is closed in $G'$.

    If $G\subset \mathbb{C}$ is open and $f:G\to \mathbb{C},f(z)=f(x+\i y)$ has
a continuous derivative, set
                 \begin{align} \partial f &=\frac{\d f}{\d z} =
                 \frac{1}{2}\left(\frac{\partial f}{\partial x}-  \i \frac
{\partial f}{\partial y}\right)&
                 \overline{\partial} f& =\frac{\d f}{\d
\overline{z}}
     =
                 \frac{1}{2}\left(\frac{\partial f}{\partial x}+  \i
\frac {\partial
  f}{\partial y}\right).\end{align}
    The operator $\partial$ is the usual derivative, $\overline{\partial}$ is called the {\em complex conjugate derivative}.
     The function $f$ is
 holomorphic if and
       only if $\overline{\partial}f =0$.
       
         In an analogous way
  one defines these derivatives for
                 Schwartz distributions. If $T$ is a distribution on $\mathbb{C}$, then $\overline{\partial} T(\varphi)=-T(\overline{\partial}\varphi)$.
 The
  function $z\mapsto 1/z$ is
                                     locally integrable and one obtains 
\begin{equation}
   \overline{\partial}(1/z)=\pi\delta(z).\end{equation}
                 Assume $f$ to be defined
and
     holomorphic for the elements $x+\i y \in G, y\ne 0$, and that $f(x\pm \i
                                      0)$
 exists and is continuous, then 
              \begin{equation}\overline{\partial}f(x+\i y)= (\i/2)(f(x+\i 0)-f(x-\i 0))\delta(y)\end{equation}
              The last equation holds for operator valued functions as well.
 To prove these assertions , we need the following variant of Gauss' theorem\\
 \begin{lem}Assume $G\subset \mathbb{C}$ an open subset , $f:G\to \mathbb{C}$ continuously differentiable and $G_0\subset G$
  an open subset with compact closure and smooth border $G_0'$, 
  then
  $$ \int_{G_0} \overline{\partial} f(z)\d^2z= -\frac{\i}{2}\int_{G_0'}f(z) \d z$$
   Here $\d z$ denotes the line element directed in the sense, that the interior of $G_0$ is on the left hand side\end{lem}

\begin{defi} Assume $R(z)$ satisfying the distribution resolvent equation, then we call
 $$ M(z)=(1/\pi)\overline{\partial}R(z)$$
the {\em spectral Schwartz distribution} of $R(z)$.
\end{defi}
\subsection{Basic properties}

If the resolvent distribution $R(z)$ is in an open set  the extension of
 the resolvent  function of an operator $A$, we call $M(z)$  a spectral distribution of $A$. Remark,
 that there may be many spectral distributions of $A$.  
  As $R(z)$ is holomorphic on an open subset, where it is continuous,   $M(z)$ is supported by the spectral set of $R(z)$.    
  
  Using the abbreviation $\r(z)=1/z$ we may the the distribution resolvent equation write in the form
  \begin{equation} R(\varphi_1)R(\varphi_2)=-R((\r*\varphi_1)\varphi_2+\varphi_1(\r*\varphi_2))\end{equation}  
  The following theorem shows, that the spectral distribution can be considered as the generalization of the family of eigenprojectors
  in finite dimensional case. The operator $M(z)$ corresponds to eigenprojecto of the eigenvalue $z$.
   We have a generalized orthogonality relation and the fact, that the product with the resolvent
  means inserting the corresponding eigenvalue.
  
  \begin{thm}\label{thm1} The spectral distribution is multiplicative. So assume  two $C_c^\infty$ functions $\varphi_1,\varphi_2$,
  then
    \begin{align*}&M(\varphi_1)M(\varphi_2)= M(\varphi_1\varphi_2)&\text{or}&\quad
             M(z_1)M(z_2)=\delta(z_1-z_2)M(z_1).\end{align*}
            Conversely
 assume a distribution $M(z)$  in $ G \subset \mathbb{C}$  obeying the equation
            $ M(z_1)M(z_2)=\delta(z_1-z_2)M(z_1)$. 
             Assume that the integral 
             \begin{align*}&
             R(z)= \int \d^2 \zeta  M(\zeta)/(z-\zeta),&&
            \int R(z) \varphi(z)\d^2 z= \int M(\zeta)\d^2 \zeta \int \varphi(z) /(z-\zeta) \d^2 z \end{align*}
            exists for $\varphi\in C_c^\infty$, then $R(z)$ fulfills the resolvent equation for distributions.
            Furthermore
            \begin{align*}&
                         M(z_1)R(z_2)=  M(z_1) /(z_2-z_1), \\&
         \int M(z_1)\varphi_1(z_1)\d ^2z_1\int R(z_2)\varphi_2(z_2) \d ^2z_2
= \int M(z_1) \varphi(z_1) \int \varphi_2(z_2)/(z_2-z_1)\d^2 z_2
            \d^2 z_1.\end{align*}
            \end{thm}
            \begin{pf}
            Use the equation
            $ \overline{\partial} \r = \pi \delta$
and hence
$ \r*\dc \varphi= \r*\dc	\delta * \varphi = \dc \delta * \r * \varphi =\dc\r*\varphi= \pi \varphi$ and
 calculate
\begin{multline*}M(\varphi_1)M(\varphi_2) = \pi^{-2} R(\dc\varphi_1 )R(\dc \varphi_2)= -\pi^{-2}R((\r*\dc \varphi_1)\dc\varphi_2
+\dc \varphi_1(\r*\dc \varphi_2))\\ 
= -\pi	R(\varphi_1 \dc \varphi_2 +( \dc \varphi_1)\varphi_2)=-\pi R(\dc(\varphi_1\varphi_2))	=M(\varphi_1 \varphi_2).\end{multline*}
In order to prove the converse assertion, observe that it means that $\r*M$ is a resolvent distribution. Now
$M(\varphi)=-1/\pi R(\dc \varphi)$ and $(\r*M)(\varphi)=-M(\r*\varphi)$. Hence
   \begin{multline*}R(\varphi_1)R(\varphi_2)=M(\r*\varphi_1)M(\r*\varphi_2)= 
   M((\r*\varphi_1)(\r*\varphi_2))\\
   =-1/\pi R(\dc((\r*\varphi_1)(\r*\varphi_2))=-R((\r*\varphi_1)\varphi_2+\varphi_1(\r*\varphi_2))\\
\end{multline*}
The last equation says 
$M(\varphi_1)R(\varphi_2) =-M(\varphi_1(\r*\varphi_2))$.
In fact
\begin{multline*}M(\varphi_1)R(\varphi_2)  = -\dfrac{1}{\pi} R(\dc \varphi_1)R(\varphi_2)=
\dfrac{1}{\pi	} R((\r*\dc \varphi_1)\varphi_2 + \dc \varphi_1 (\r*\varphi_2))\\=
\dfrac{1}{\pi	} R(\pi \varphi_1\varphi_2 + \dc \varphi_1 (\r*\varphi_2))
=\dfrac{1}{\pi}R(\dc (\varphi_1(\r*\varphi_2))=-M(\varphi_1(\r*\varphi_2)) .\end{multline*}
\end{pf}
The operator $M(z)$ is a generalized eigen-projector. In fact if $R(z)=1/(z-A)$ then\\
$M(z_1)(1/(z_2-A) =  M(z_1) /(z_2-z_1).$ The relation $M(z_1)M(z_2)=\delta(z_1-z_2)M(z_1)$ is a generalized orthonormality relation.
    \begin{rem} \label{rem1}If $A \in L(V)$ is an operator and
     if  the distribution  $R(z)$ fulfills  the equation $AR(z)=-1+zR(z)$ in the sense of distributions, then it fulfills the
     resolvent  equation (\ref{eq4})
            $ R(z_1)-R(z_2)=(z_2-z_1)R(z_1)R(z_2).$
            and it is not necessarily the resolvent distribution of $A$ . If it is a resolvent distribution , by differentiation one obtains in this case
            $ AM(z)= zM(z)$,
            as $ \overline{\partial}z=0$.So $M(z)$ is an eigenprojector  for the eigenvalue $z$.\end{rem}

One proves by a partition of unity
\begin{prop}Assume a   bounded set $G \subset \mathbb{C}$ and a family $G_i,i=1,\cdots,l$ of open, pairwise disjoint subsets of $G$,
furthermore an open subset $G_0 \supset  (G \setminus \bigcup_1^l G_i)$ and assume a distribution $M$ on $G$, whose restriction
to $G_0$ vanishes and whose restriction to $ G_i$ is multiplicative for $i=1,\cdots,l$.
 Then $M$ is multiplicative on $G$ if and only if  $M(\varphi)M(\psi)=0$ ,  under the condition that
 $\varphi$ has its support in $  G_i$ and $\psi$ has
  its support in $ G_j$ with
$i,j\in [1,l],i\ne j$.
 
 \end{prop}
 
\begin{prop}\label{prop2} 
 Assume an open bounded  set $G \subset \mathbb{C}$ and a family $G_i,i=1,\cdots,l$ of open, pairwise disjoint subsets of $G$.
  Assume in each subset $G_i$ a subset  $G_{i,0}$ of
 Lebesgue measure $ 0$ and a resolvent function $R(z)$ on $G \setminus \bigcup G_{i,0}$ and  for each $i$ a sequence of subsets
   $G_i \supset G_{i,n}\downarrow G_{i,0}$, such that 
   $$ R_i(\varphi)= \lim_{n \to \infty}\int _{G_{i}\setminus G_{i,n}} R(z) \varphi(z)  \d^2  z$$
     for $\varphi \in \mathcal{D}(G_i)$ exists and defines a resolvent distribution on $G_i$.
Then the distribution $R$ defined by the  the function
     $R(z)$ on $G \setminus \bigcup G_{i,0}$ and by the distributions $R_i$ on $G_i$ fulfills the resolvent distribution equation.
     \end{prop}
     
  \begin{pf}   
 Call $M_i= 1/ \pi \dc R_i$. We have to show, that 	$M_i(\varphi_1)M_j(\varphi_2)=0$ , if the support of $\varphi_1$ is in
 $G_i$ and the support of $\varphi_2$ is in $G_j$ for $ i \ne j$. We have using (\ref{eq5}) 
 $$ \int_{G_{i}\setminus G_{i,m}}\d^2 z_1 \int_{G_{j}\setminus G_{j,n}}\d^2 z_2\varphi_1(z_1)\varphi_2(z_2)R(z_1)R(z_2) = I+II$$
 with
 $$ I= -\int_{G_i\setminus G_{i,m}}\d^2 z_1 \int_{G_j\setminus G_{j,n}}\d^2 z_2\varphi_1(z_1)\varphi_2(z_2)R(z_1)\frac{1}{z_1-z_2}=
 -\int_{G_i\setminus G_{i,m}}\d^2 z R(z) \psi_n(z)$$	
 The sequence
 $$\varphi_1(z_1) \int_{G_j \setminus G_{j,n}}\d^2 z_2\varphi_2(z_2)/(z_2-z_1)=\psi_n(z_1)
 \to \varphi_1(z_1) \int_{G_j }\d^2 z_2\varphi_2(z_2)/(z_2-z_1)=\psi(z_1)$$
in the sense of $\mathcal{	D}(G_i)$.
 We have 
  $ I \to -R_i(\psi)=-R(\psi)= -R(\varphi_1(\r*\varphi_2)).$
 Similarly $II \to -R(\varphi_2(\r*\varphi_1))$.
  By the proof of  theorem \ref{thm1}, one sees that $M_i(\varphi_1)M_j(\varphi_2) = M(\varphi_1\varphi_2)
 =0$.
 \end{pf}

 
 \begin{prop}\label{prop3}
  If the resolvent $R(z)$  has a  pole of order $n$  in a point $z_0$, then $R(z)$ is holomorphic for $z$ in a neighborhood of $z_0, z \ne z_0$
 and behaves like $ b_0( z-z_0)^{-1} + \cdots b_{n-1} (z-z_0)^{-n}$  near $z_0$. We have $b_0=p$ and $b_j= a^j$, where $p^2=p$ and
 $a^n=0$ and $ap=pa$. Define
 $$ R(z) = \sum_{k=0}^{n-1} b_k\frac{\mathcal{P}}{(z-z_0)^{k+1}},$$
 where the distribution \cite{Schwartz1}
 $$ \int \d^2 z \frac{\mathcal{P}}{z^k}\varphi(z)=\lim_{\varepsilon \downarrow 0} \int_{|z|> \varepsilon} \frac{1}{z^k}\varphi(z)\d^2 z
 = (-1)^k \int \frac{1}{z}\partial ^k \varphi(z)\d^2 z.$$
 The distribution $R(z)$ fulfills the resolvent distribution equation. The spectral distribution is
 $$ M(z)=  \sum_{k=0}^{n-1} \frac{(-1)^k}{k!} p a^k \partial^k \delta(z-z_0)$$
 \end{prop}
 The proof follows directly from Leibniz's formula. 
 \begin{prop}\label{prop4}
 Assume a compact set $K$ in the real line and an open neighborhood $G=G_0 \times I$ of $K$, where  $G_0$ is an open set containing
 $K$ and $I$ is an open interval containing $ 0$
,   and a resolvent function $R(z)$
 holomorphic in $G\setminus K $. Define the function $R_u$ on $I$ by $R_u(x)= R(x+\i u)$ and assume, that for any 
 $C^\infty$-function $\varphi$ with support in $G_0$ the limits of $R_u(\varphi)= \int \d x R_u(x) \varphi(x)$ for $ u \to 0+$
 and for $u \to 0-$ exist. If $\varphi$ is a $C^\infty_c$- function on $G$, define $\varphi_u(x)=\varphi(x+ \i u)$. Then
 $ u \mapsto R_u(\varphi_u)$ is integrable and defines
   the distribution  $R(\varphi)$ on $G$ by
 $$ R(\varphi)= \int \d u  R_u(\varphi_u) = \int \d u \int \d x R(x+ \i u) \varphi(x+\i u).$$
 The integral defines a resolvent distribution on $G$.
 \end{prop}
 
\begin{pf} 
 Define for $m=0,1,\cdots$ and for any test function $\varphi$ on $G_0$ the norm
 $$ \| \varphi \| _m = \max{\{|\partial_x^k\varphi(x)|:	x\in G_0,k=1,\cdots,m\}}.$$
 For a distribution $T$ on $ G_0$ define accordingly $ N_m(T)= \sup \{\|T(\varphi)\|:\| \varphi\| _m\le 1\}$
 Denote by $R_{0\pm}$ the limites of $R_u$.
 As for distributions weak convergence implies strong convergence , there exists an $m$ such that $N_m(R_u-R_{0+})\to 0$
 and $ N_m(R_u-R_{0-})\to 0 $  for $u \to 0+$ or $u \to 0-$ resp.. The function  $ u \mapsto R_u(\varphi_u)$ is continuous
 for $u \ne 0$ and has right and left limits at $u=0$. Hence it is integrable. Define a $C^\infty$ function $ \alpha $ on the  real line with
 $\alpha \ge 0, \alpha (x)=0 \text{ for } |x| \ge 1, \alpha(x) =1 \text{ for} |x| \le 1/2$. Put $ \alpha_\varepsilon(u)= \alpha(u/\varepsilon)$.
 Then for $z_j= x_j+ \i u_j, j=1,2$
 $$ R(\varphi_1)R(\varphi_2)=\iint \d u_1 \d u_2 R_{u_1}(\varphi_{1,u_1})R_{u_2}(\varphi_{2,u_2})
 = \lim _{\varepsilon\to 0}\iint \d u_1 \d u_2 R_{u_1}(\varphi_{1,u_1})R_{u_2}(\varphi_{2,u_2})(1-\alpha_\varepsilon(u_1-u_2)) $$
 On the other hand
 
$$ \iint \d u_1 \d u_2 R_{u_1}(\varphi_{1,u_1})R_{u_2}(\varphi_{2,u_2})(1-\alpha_\varepsilon(u_1-u_2))=
\iint\d^2 z_1 \d ^2z_2 R(z_1)R(z_2) \varphi_1(z_1)\varphi_2(z_2)(1-\alpha_\varepsilon(u_1-u_2)).$$
In the second term the quantities $R(z_1),R(z_2)$ are usual resolvent functions. We apply the resolvent equation (\ref{eq5}) and split
the integral into two terms $I_\varepsilon$ and  $II_\varepsilon$. 
$$I _\varepsilon= -\iint \d^2 z_1 d^2 z_2 R(z_1)1/(z_1-z_2) \varphi_1(z_1)\varphi_2(z_2)\alpha_\varepsilon(u_1-u_2)
= -\int \d u R(\psi_{u,\varepsilon})$$
with 
$$ \psi_{u,\varepsilon}(x) = \varphi_1(x + \i u) \int \d x_2 \d u_2  \varphi_2(x_2+ \i u_2)\frac{1}{x+ \i u -x_2-\i u_2}
(1- \alpha_\varepsilon(u-u_2)).$$
Put
$$ \psi_u(x)=\varphi_1(x + \i u) \int \d x_2 \d u_2  \varphi_2(x_2+ \i u_2)\frac{1}{x+ \i u -x_2-\i u_2}.$$
We have to estimate
$$ \partial_x^k(\psi_u(x)-\psi_{u,\varepsilon}(x))=\sum_{j=1}^k \binom{k}{j}
\partial^{k-j}_x\varphi_1(x + \i u) \int \d x_2 \d u_2 \partial^j_x \varphi_2(x_2+ \i u_2)\frac{1}{x+ \i u -x_2-\i u_2}
\alpha_\varepsilon(u-u_2)
$$
The integral can be estimated with some constants $C,K$
$$|\int \d x_2 \d u_2 \partial^j_x \varphi_2(x+\i u-x_2- \i u_2)\frac{1}{x_2+\i u_2}
\alpha_\varepsilon(u_2)|\le
 C \iint_{|x'| \le K, |u'|\le \varepsilon}\d x' \d u' \frac{1}{\sqrt{x'^2+u'^2}}= O (\varepsilon| \ln \varepsilon|)$$ 

 Hence 
 $ \|\psi_u-\psi_{u,\varepsilon}\| _m=O(\varepsilon|\ln \varepsilon|)$ uniformly in $u$ and
 $$ I_\varepsilon \to I = -\iint \d^2 z_1 d^2 z_2 R(z_-1)1/(z_1-z_2) \varphi_1(z_1)\varphi_2(z_2)$$
     By interchanging the roles of 1 and 2, one obtains the corresponding result for $II_\varepsilon$.
     \end{pf}
   \begin{thm}\label{thm2}
   If $R$ is the resolvent of a bounded operator and $M$ is a resolvent distribution extending $R$ defined
   on the whole complex plane , then the support of $M$ is compact  and $M(1)$ is defined , where here $1$ is the constant function $1$,
   and
   $$M(1)=\int \d^2 z M(z) =1_{L(V)}$$
   In this case we call $M$ {\em complete}.
   \end{thm}
   \begin{pf}Assume, that the support of $M$ is contained in the circle of radius $r$. Assume a test 
    function $\varphi$	constant $1$ on this
   circle, then 
   \begin{multline*}
   $$ M(1)=M(\varphi)= -\frac{1}{\pi}\int  \d^2 z R(z) \dc \varphi(z)= -\frac{1}{\pi}\int _{|z|>r} \d^2 z R(z) \dc \varphi(z)=\\
    -\frac{1}{\pi}\int_{|z|>r} \d^2 z \dc( R(z)  \varphi(z))=
   \frac{1}{2\pi\i}\int_\Gamma \d z R(z)= 1_{L(V)},\end{multline*}
   where $\Gamma$ is the circle of radius $r$ run in the anti-clockwise sense.
    Now the residuum  in infinity of $R(z)=1/(z-A)$ equals $1_{L(V)}$.
   \end{pf}

\subsection{Examples} \begin{enumerate} \item    
           Finite dimensional matrix.
      By Jordan's normal form one obtains that
               $$ M(z)=\sum _i p_ i \sum_k   (1/k!)(-1)^k a_i^k \partial^k \delta(z- \lambda_i),$$
               where the $\lambda_i$ are the eigenvalues,the $p_i$ are the eigenprojectors,
               $p_ip_j= \delta_{ij}$ and the $a_i$ are nilpotent and $a_ip_j= \delta_{ij}a_i$.   Instead of $\partial$ we could have 
               chosen any other linear combination
                $D$ of $\partial_x$ and $\partial_y$, such that $Dz=1$. This an example, that there are many resolvent distributions extending
                a resolvent function.
               
               We cite another example affirming remark \ref{rem1}.  Assume $A=0$, then 
               $zR(z)=1$, so $R(z)=1/z-\pi	C \delta(z)$, where $C$ is an arbitrary matrix.  Then  $M(z)= \delta(z)-C \overline{\partial}\delta(z) $ and
               $M(\varphi)M( \psi)=M(\varphi \psi)$  if and only if $C^2=0$, this might not be the case.

               \item Assume that $V$ is a Hilbert space and that $U$ is a unitary operator.
                  $$ \int M(z) \varphi(z) \d^2 z = \sum_{l= -\infty}^\infty U^l \frac{1}{2 \pi}\int \e^{-\i \vartheta l} \varphi (\e ^{\i \vartheta })
        \d \vartheta. $$
         \item Assume $V$ to be a  Hilbert space and $A$ to be a  self adjoint  operator.  Let  $E(x), x \in \mathbb{R}$ be the
        spectral family of $A$,then$$ M(x + \i y) = E'(x)\delta(y), $$
              where the derivative is in the sense of distributions.
              \end{enumerate}

        \section{Discussion of the Multiplication Operator and some Rank One  Perturbations}
  \subsection{ $C^\infty$ Multiplication Operator}
        Assume an open  bounded set
         $ G \subset \mathbb{R}^n $  and a $C^\infty$  bounded function $P:G\to \mathbb{R}$ ,
         which is bonded with all its derivatives and with $ \nabla P(y) \ne 0$ for all $ y \in G$.
Consider   $ V= L^2(G)$ and the operator $\Omega: V\to V$ with $ \Omega f(y)=P(y)f(y)$ . The operator $\Omega$ is bounded
and its resolvent is
 $R(z)= 1/(z-\Omega)$  for 
$\Im z \ne 0$. 
      For any test function on the complex plane the integral
      \begin{align*}&
         R(\varphi)= \int d^2 z \varphi(z)/(z-\Omega),&&
           ( R(\varphi)f)(y)= \int d^2 z \varphi(z)/(z-P(y))f(y)\end{align*}
        exists and we define the resolvent distribution in that way. A short calculation shows, that the resolvent distribution equation is fulfilled.
We find for the spectral distribution
         $$ M(\varphi)f(y)= \varphi(P(y))f(y).$$
               The open mapping theorem ensures that the set
         $ P(G)=\{ x\in \mathbb{R}: \exists y\in G, P(y)=x\}$ is open. The set
         $$ S(x)= \{ y \in G:P(y)=x\}$$ is a $(n-1)$-dimensional $C^\infty$-submanifold of $\mathbb{R}^n$. With the help 
         of the theorem of implicit functions one establishes the lemma
\begin{lem}\label{lem2}
         Assume a point $y_0\in G$ and set $x_0=P(y_0)$. The exists an open neighborhood $U_1 \subset \mathbb{R}$ of $x_0$
          and an open subset $U_2 \subset\mathbb{R}^{n-1}$ and an  injective $C^\infty$ mapping  $\Psi: U_1\times U_2 \to G$, such that
         $$P( \Psi(x,u))=x$$ 
         and $\Psi: U_1\times U_2\mapsto \Psi(U_1\times U_2)$ is a diffeomorphism. 
        Then $\Psi(U_1\times U_2) $ is an open neighborhood of $y_0$ and
         for fixed $x$ the mapping $\Psi_x: u \in U_2 \mapsto \Psi(x,u) $ is a $C^\infty$-chart of $S(x)$. 
         If
               $$ J(\Psi)= | \det ( \partial \Psi/\partial x, (\partial \Psi/\partial u_i)_{i=1,\cdots,n-1}))|$$
       is the absolute value of the Jacobi's determinant and
       $$ \Gamma (\Psi_x)=( \det (\partial \Psi_x/ \partial u_i.\partial \Psi_x/ \partial u_j)_{i,j=1,\cdots,n-1})^{1/2}$$
       is the square root of Gram's determinant and $\d \sigma_x(y)$ is the euclidean surface element  on $S(x)$ for $y\in S(x)$,
        where the point denotes the scalar product. 
        Then 
        $ J(\Psi)(x,u)= \Gamma(\Psi_x)(x,u)/| \nabla P(x,u)|$
        and
        $ \d y = \d x \ \sigma_x(y)/|\nabla P(y)|$
     \end{lem}

              We introduce 
       $ \d \tau_x(y)= \d \sigma_x (y)/|\nabla P(y)|$.
       So 
       \begin{equation}\d y = \d x \d \sigma_x(y)/|\nabla P(y)|= \d x \d \tau_x(y)\end{equation}
         
       As any function of compact support on $G$ can be represented as a finite sum of functions which have their support in a chart
       like in the preceding lemma, we obtain the proposition
       
       \begin{cor}\label{cor1} If $f$ is a continuous function of compact support on $G$, then
        $$   \int f(y) \d y =  \int \d x \int_{y \in S(x)}\d  \tau_x(y) $$
        \end{cor}   
        \begin{lem} Assume a function $f\in \mathcal{D} (G)$, then
       $$ \omega \in P(G) \mapsto  \int_{y \in S(\omega)} f(y) \d \tau_x (y) $$
       is in $\mathcal{D}(P(G))$where $P(G)$ is the image of $G$ in $\mathbb{R}$ \end{lem}
       \begin{pf}
       Assume at first, that the support of $f$ is contained in an open set $\Psi(U_1\times U_2)$  like in lemma \ref{lem2} 
       . Using again this lemma 
       we have
       $$  \int_{y \in S(\omega)} f(y) \d \tau (y)= \int \d u f(\Psi(\omega, u) J(\omega,u) $$
       and this is surely $C^\infty$. A partition of unity finishes the proof. \end{pf}
       
       If $T$ is a distribution on the open set $P(G)$, then consider $T(P(y))$. If $f$ has its support in $U_1 \times U_2$, then as 
       distributions transform like functions
      $$ \int T(P(y)) f(y) \d y= \iint \d x \d u T(P(\Psi(x,u)) f(\Psi(x,u)) J(\Psi)(x,u)= 
       \int \d x T(x) \int \d u f(\Psi(x,u))\Gamma(\Psi_x)(u)$$
       Hence we define for test function $ f \in \mathcal{D}(G)$ 
       \begin{equation}T(\Omega)(f)=\int T(P(y)) f(y) \d y= \int \d x T(x) \int_{y \in S(x)} f(y) \d \tau_x (y).\end{equation}
      Especially
      \begin{equation} \delta(x-\Omega)(f) = \int_{y \in S(x)}  f(y) \d \tau_x (y)\end{equation}
       and $M(z)=\delta(x-\Omega)\delta(y)$. This measure on $\mathbb{R}^n$ has been treated by Gelfand-Schilow  under the 
       name $\delta(x-P)$\cite{Gelfand1}.
        We cite the definition \cite{Gelfand2}
         \begin{defi} Assume a vector space  $V$ and and a linear mapping $A:V \to V$. A linear functional $F:V \to \mathbb{C}$ is called
         a generalized eigenvector for the eigenvalue $x$ if $F(Af)=xF(f)$ for all $f \in V.$.
         \end{defi}
  Obviously  the operator
         $\Omega$ leaves $ \mathcal{D}(G)$ invariant. We define the left generalized ket-vector $\langle \delta_y|$ as the kernel
         $\langle \delta_y| (w)=\delta(y-w)$
         applied to $f$ in the following way
         $$\langle  \delta_y|f \rangle = \int \delta(y-w) f(w) \d w =f(y)$$
         hence
         $$ \langle \delta_y| \Omega f\rangle =P(y) f(y)= P(y)\langle \delta_y |f \rangle $$
         and $ \langle \delta_y|$ is  a generalized left eigenvector for the eigenvalue $P(y)$.
      . Similar we define the generalized  right eigenvector $| \delta_y \rangle$
         by $\langle f|\delta_y \rangle =  \overline{f}(y)$.
         For $f \in L^2(G)$ the bra-vector $\langle f|$ is given by the functional $ g \mapsto \langle f|g =\langle f| g \rangle = 
         \int \overline{f}(y)g(y) \d y$. Similar relations hold for the ket-vector $|f \rangle$.  Following this idea we define
         $ \langle \delta_y| \delta_{y'} \rangle = \int \d w \delta(y-w)\delta(y'-w) $ and finally
         $$ \langle \delta_y| \delta_{y'} \rangle= \delta(y-y').$$	
         This is the {\em orthogonality relation} for the generalized eigenvectors.
        The relation 
         $$ \langle f|g \rangle =\iint \d y \d y' \overline{f}(y)g(y') \langle \delta_y| \delta_{y'}\rangle.$$
         states the {\em completeness} of the eigenvectors.
         Finally
         $$ \delta(x-\Omega) = \int _{y \in S(x)} |\delta_y \rangle \langle \delta_y | \d \tau_x(y).$$

       \subsection{Perturbations of the Multiplication Operator 1}
       We perturb the multiplication operator $\Omega$
       and define the operator
       \begin{equation} H=\Omega+ |g \rangle \langle h |\end{equation}
       with $g,h \in \mathcal{D}(G)$. The resolvent can be easily calculated   
       \begin{align}\label{eq14}&                               
        R(z)=R_\Omega(z)+ R_\Omega(z)|g \rangle \langle h| R_\Omega(z)/C(z),&&
        C(z) = 1- \langle h| R_\Omega (z)|g \rangle, \end{align}
       where $R_\Omega(z)=1/(z-\Omega)$ is the resolvent of $\Omega$ .
       The ket-vector $  R_\Omega(z)|g \rangle$ is the functional $ f\in \mathcal{D}(G) \mapsto \langle f|  R_\Omega(z)|g \rangle$ and similarly
       defined is the bra-vector $ \langle h| R_\Omega(z)$.
       
        Recall the formula
       $$ \frac{1}{x\pm \i \varepsilon}\to \frac{1}{x\pm \i 0} = \frac{\mathcal{P}}{x} \mp \i \pi\delta(x) \text { for  } \varepsilon \downarrow 0,$$
      where $\mathcal{P}$	denotes the principal value and obtain the lemma
       
             \begin{lem}\label{lem4}
Assume $f_1,f_2 \in \mathcal{D}(G)$, then for $z \to x \pm \i 0$ uniformly in $x$
$$\int \d y \frac{\overline{f}(y)g(y)}{z-P(y)}= \langle f_1| \frac{1}{z-\Omega}|f_2 \rangle \to 
\langle f_1|\frac{1}{x\pm \i 0-\Omega}|f_2 \rangle = 
\langle f_1| \frac {\mathcal{P}}{x-\Omega}
 |f_2 \rangle 
\mp \i \pi \langle f_1| \delta(x- \Omega) |f_2 \rangle $$
with
\begin{align*}\langle f_1| \frac {\mathcal{P}}{x-\Omega} |f_2 \rangle=&
\int \d \omega \frac {\mathcal{P}}{x- \omega}\int_{y \in S(\omega)}\overline{f} _1(y)f_2(y) \d \tau (y),&
\langle f_1| \delta(x- \Omega) |f_2 \rangle=& \int_{y \in S(x)}\overline{ f}_1(y)f_2(y) \d \tau (y)
\end{align*}
\end{lem}
We use Dirac's notation $\langle f_1|Af_2 \rangle=\langle f_1|A|f_2 \rangle$, which is in many cases convenient.
We calculate for $ z \to x \pm \i 0 $
$$ C(z)\to  
1-  \langle h| \frac {\mathcal{P}}{x-\Omega}| g \rangle 
\pm \i \pi \langle h| \delta(x- \Omega) |g \rangle = C_1(x) \pm \i \pi C_2(x) $$

        So there are jumps of $C(z)$
          on the real axis contained in $P(K)$, if the supports of $g.h$ are contained in the compact set $K \subset G$. 
           If $ z \notin P(K)$ and $C(z) \ne 0$, then $R(z)$ exists as function and
        $z$ belongs to the resolvent set. We assume, that the jumps on $P(K)$ are 
       are isolated,  more precisely that there exists an open neighborhood  $G_0 \subset \mathbb{C}$ of $P(K)$, such that 
       $C(z)$ is $\ne 0$ and holomorphic in $G_0 \setminus  P(K)$  and $C(x \pm \i 0)\ne 0$ for $x \in G_0 \cap \mathbb{R}$.
       We may assume, that $G_0=G_1\times I$, where $G_1$ is an open subset of the real line containing $P(K)$ and $I$ is an open interval
       containing 0.
       
       \begin{lem}\label{lem5}    Assume a  $ C^\infty$ test function $\varphi$ with support in $G_1$
        and put $R_u(x)=R(x+ \i u)$
       then for $\varphi \in   \mathcal{D}(G_1)$	the operator $  \int \d x R(x+ \i u) \varphi(x)= R_u(\varphi)$ on $L^2(G)$
       converges for $u \to 0+$ or $u \to 0-$
       in operator norm to an operator called $R_{\pm 0}(\varphi)$.
       \end{lem}
    
\begin{pf}
       Rewrite equation (\ref{eq14})
      $$ R(z)= \frac{1}{z- \Omega} + \iint \d w_1 \d w _2 \frac{1}{z-w_1}\frac{1}{z-w_2} \delta(w_1-\Omega)|g \rangle \langle h |
       \delta(w_2-\Omega)/C(z)$$
     $$  = \frac{1}{z- \Omega} + \iint \d w_1 \d w _2 \frac{1}{w_2-w_1}(\frac{1}{z-w_1}-\frac{1}{z-w_2}) \delta(w_1-\Omega)|g \rangle \langle h |
       \delta(w_2-\Omega)/C(z)$$

       Then for $ u \ne 0$
$$\int \d x R(x+ \i u) \varphi(x)= R_u(\varphi)= \psi_u(\Omega)+ \iint \d w_1 \d w_2 \chi_u(w_1,w_
2) \delta(w_1-\Omega)|g \rangle \langle h |
       \delta(w_2-\Omega)$$
       The function $\psi_u(w)= \int \d x  \varphi(x)/(x + \i u -w)$ is in $\mathcal{D}(G_1)$ and
         converges in this sense to $ \psi_{+0}$ or $\psi_{-0}$.
      $$ \chi_u(w_1,w_2)= \int \d x \varphi(x)\frac{1}{w_2-w_1}
     ( \frac{1}{x + \i u-w_1}-\frac{1}{x+ \i u-w_2})/C(x + \i u)$$ 
     $$ =- \int \d x \frac{ \varphi(x)}{C(x+ \i u)}\int _0^1 \d t \frac{\d}{\d x} \frac{1}{x+ \i u -w_1- t(w_2-w_1)}$$
     $$ = \int \d x \frac{\d}{\d x}(\frac{\varphi(x)}{C(x+ \i u)})\int _0^1 \d t  \frac{1}{x+ \i u -w_1- t(w_2-w_1)}
     = \iint \d x \d t \frac{ \omega_u(x+(1-t)w_1+tw_2))}{x+ \i u} $$
     $$\omega_u(x)=  \frac{\d}{\d x}(\frac{\varphi(x)}{C(x+ \i u)})$$
     The function  $x \mapsto \omega_u(x)$ is in $\mathcal{D}(G_1)$ and converges for $ u \to \pm 0$ in this sense , as the
     function $x\mapsto C(x +\i u)$ is in $C^\infty$ and converges 
    to $x \mapsto C(x\pm 0)$ for $x\to +0$ or $x \to-0$
     uniformly , the analogue holds for all derivatives. Hence    
     we obtain
      $$ \chi_{\pm 0}(w_1,w_2)= \iint \d x \d t \frac{ \omega_ {\pm 0}(x+(1-t)w_1+tw_2))}{x\pm  \i 0}. $$
      Assume $f_1,f_2\in L^2(G)$, then for $\vartheta = \pm 0$
      \begin{multline*}
       |\iint \d w_1 \d w _2( \chi_u(w_1,w_2)- \chi_\vartheta (w_1,w_2)) \langle f_1|\delta(w_1-\Omega)|g \rangle
      \langle h| \delta(w_2-\Omega)| f_2 \rangle|\\
      \le \max{\{ | \chi_u(w_1,w_2)- \chi_\vartheta (w_1,w_2))|:w_1,w_2 \in P(K)\}} \iint \d w_1 \d w_2 
      | \langle f_1|\delta(w_1-\Omega)|g \rangle
      \langle h| \delta(w_2-\Omega)| f_2 \rangle|.\end{multline*}
      Observe, using corollary \ref{cor1} , that e.g. $\int \d w |\langle f| \delta(w-\Omega)|g \rangle| \le N(f)N(g)$,
        where $N(.)$ is the $L^2$-norm,
      and conclude from there that
      the last expression converges to 0 in operator norm,
      \end{pf}
      A consequence of proposition \ref{prop4} and lemma \ref{lem5}
      \begin{prop}
           Define for a  $ C^\infty$ test function $\varphi$ with support in $G_0$ the operator valued distribution
       $$R(\varphi) = 
       \lim _{\varepsilon \downarrow 0}  \iint _{|u| > \varepsilon} \d x \d u R(x,u)\varphi(x,u)=
       \lim _{\varepsilon \downarrow 0}  \int _{|u| > \varepsilon}  \d u R_u(\varphi_u)$$
       with $\varphi_u(x)=\varphi(x+\i u)$.
       Then $R(\varphi)$ extends the function $R(z)$ on $G_0 \setminus P(K)$  to a distribution on $G_0$, which fulfills 
       the resolvent distribution equation.
       We have        
 \begin{align*} M(x+\i y)=&\mu(x) \delta(y),&\mu(x)= &\frac{1}{2 \pi \i}(R(x-\i 0)-R(x+ \i 0). \end{align*}
       \end{prop}
       We consider for $u  \ne 0$ the sesquilinear form
$$ f_1,f_2 \in \mathcal{D}(G) \to  \mathcal{B}(R_u(x))(f_1,f_2)= \langle f_1| R_u(x)| f_2 \rangle. $$
By lemma \ref{lem4} the bracket converges uniformly in $x$ for $ x \to \pm 0$ to a sesquilinear form $\mathcal{B}(R_{\pm 0})(x)$. We obtain
$$ \langle f_1| R_{\pm 0}(\varphi)| f_2 \rangle = \int \d x \varphi(x) \mathcal{B}(R_{\pm 0})(x)(f_1,f_2),$$
and define
 $$     \langle f_1 | \kappa(x)|f_2 \rangle= 	 \mathcal{B}(\mu(x))(f_1,f_2)=
 \frac{1}{2 \pi \i}(\mathcal{B}(R_{-0})(x)(f_1,f_2)-\mathcal{B}(R_{+0}(x)(f_1,f_2)).$$
 Hence
 $$ \langle f_1| \mu(\varphi)| f_2 \rangle = \int \d x \varphi(x)\langle f_1 | \kappa(x)|f_2 \rangle $$
 and $\kappa(x)$ appears as the restriction of the sesquilinear form $f_1,f_2 \in L^2(G) \mapsto  \langle f_1| \mu(x)| f_2 \rangle$, which 
 is a distribution, to $f_1,f_2 \in \mathcal{D}(G)$,  where $ x \mapsto  \langle f_1 | \kappa(x)|f_2 \rangle$ is a continuous  function. 
 
From the formula
$$ \mathcal{B}(R_{\pm 0}(x))=  \frac{\mathcal{P}}{x-\Omega} \mp \i \pi \delta(x-\Omega)+
\frac{1}{C_1\pm\i \pi C_2}(\frac{\mathcal{P}}{x-\Omega}\mp\i \pi \delta(x-\Omega))|g \rangle
\langle h|( \frac{\mathcal{P}}{x-\Omega}-\mp\i \pi \delta(x-\Omega)). $$
 one establishes by straight forward calculations

\begin{prop}
Recall
\begin{align*} C_1(x)=&1-  \langle h| \frac {\mathcal{P}}{x-\Omega}| g \rangle 
& C_2(x)=& \langle h| \delta(x- \Omega) |g \rangle 
\end{align*}
Assume, that $g,h$ are real and  that 
$$ |C(x \pm \i 0)|^2= |C_1(x)  \pm \i \pi C_2(x)|^2 = C_1(x)^2+ \pi^2 C_2(x)^2 \ne 0$$
for $x \in  G_0\cap \mathbb{R}$. 
 We define the following bra- and ket-vectors  as functionals over $\mathcal{D}(G)$.
           \begin{align*}
            A=&A(x)= \frac{\mathcal{P}}{x-\Omega} |g \rangle &
        A'=&A'(x)= \langle h| \frac{\mathcal{P}}{x-\Omega}  \\
        B=&B(x)=\delta(x-\Omega)|g \rangle &
               B'=&B'(x)=\langle h|\delta(x-\Omega) \end{align*}     
                For $x \in P(K)$ one obtains
       $$\kappa(x)	 = \delta(x-\Omega)+ \frac{1}{C_1^2+\pi^2 C_2^2}
       (AC_2A'+AC_1B'+BC_1A'- \pi^2 BC_2B')$$
                If $C_2(x) \ne 0$, we may write
       $$\kappa(x)= \delta(x-\Omega)-\frac{BB'}{C_2}
+ \frac{1}{(C_1^2+\pi^2 C_2^2)C_2}
       (AC_2^2A'+AC_1C_2B'+BC_1C_2A'+ C_1^2BB')$$ $$
       = \delta(x-\Omega)-\frac{BB'}{C_2}+ \frac{1}{(C_1^2+\pi^2 C_2^2)C_2}
       ( C_1B+AC_2)(C_1B'+A'C_2)
        = p(x)+|\alpha(x) \rangle \langle \alpha'(x)|$$ 
        with
        $$ p(x)= \delta(x-\Omega)-\frac{\delta(x-\Omega)|g \rangle \langle h|\delta(x-\Omega)}
        {\langle h|\delta(x-\Omega)|g \rangle} $$
      $$|\alpha(x) \rangle = N(x) (C_1(x) \delta(x-\Omega)|g \rangle+ C_2(x)\frac{\mathcal{P}}{x-\Omega} |g \rangle)$$ 
    $$ \langle \alpha'(x)|= N(x)(C_1(x)\langle h| \delta(x-\Omega)+C_2(x)\langle h|\frac{\mathcal{P}}{x-\Omega}$$
 and $N(x)^2 =1/((C_1^2+\pi^2 C_2^2)C_2)$.
 \end{prop}
       We discuss the case $C_2(x)=0$.  If $h=g$, then $C_2(x)= \int_{S(x)} \d \tau(y) g(y)^2 =0$ implies , that $g$ on $S(x)$
       vanishes and hence $B=0,B'=0$. In the general case, remark that $B=0$ or $B'=0$ imply that $C_2=0$.  Let us assume, that
       $g,h$ are in such form ,that $C_2(x)=0$ implies $B(x)=0,B'(x)=0$. Then
       \begin{equation} \kappa(x)=\begin{cases}p(x)+|\alpha(x) \rangle \langle \alpha'(x)|& \text{ for } C_2(x) \ne 0\\
       \delta(x-\Omega) & \text{ for } C_2(x)=0\end{cases}\end{equation}     
         If $g=h$, then $\mu(x)$ is positive definite.

         \begin{prop}
         The operator $H$ maps $\mathcal{D}(G) \to \mathcal{ D}(G)$. So we may formulate
         \begin{align*}
         & H \kappa(x)= \kappa(x) H =x \kappa(x)&&
         Hp(x)=p(x)H =xp(x)\\&
         H|\alpha(x) \rangle = x |\alpha(x) \rangle&&
          \langle \alpha'(x)|H =x \langle \alpha'(x)|\end{align*}
         These  equations have to be understood as equations for sesquilinear forms bracketed between functions in $\mathcal{D}(G)$ or
         as functionals on $\mathcal{D}(G)$ . So $\kappa(x)$ and $p(x)$   are generalized eigenprojectors and $|\alpha(x) \rangle$ is a generalized
         right eigenvector and $\langle \alpha'(x)|$ is a generalized left eigenvector,  all for the eigenvalue $x$. \end{prop}
         The proof is done by straight forward computation.
         \begin{prop}\label{prop8}
         We have the orthogonality relations
         \begin{align*}
         & \kappa(x)\kappa(x')= \delta(x-x')\kappa(x)&&
        p(x)p(x')= \delta(x-x')p(x)\\&
        p(x)|\alpha(x') \rangle =0, \;\langle \alpha'(x)| p(x')=0&&
        \langle \alpha'(x)| \alpha(x')\rangle= \delta(x-x')\end{align*}
       These relations have to be understood in the sense of distributions, e.g. the first one signifies
       $$(\int \d x_1 \varphi_1(x_1) \langle f_1| \kappa(x_1)) (\int \d x_2 \varphi_2(x_2)  \kappa(x_2)| f_2 \rangle
         =\int \d x \varphi_1(x)\varphi_2(x)\langle f_1| \kappa(x)| f_2 \rangle $$
         for $f_1,f_2 \in \mathcal{D}(G)$ and $\varphi_1,\varphi_2 \in \mathcal{D}(P(G))$ and we show, that the expressions make sense.
       
       \end{prop}
       \begin{pf}
          Assume $f_1 \in \mathcal{D}(G)$, by the mapping $f_2 \in \mathcal{D}(G)\to \langle f_1| \kappa(x)| f_2 \rangle$ we define a   distribution 
           called  $ \langle f_1 | \kappa(x)$.  For a test function $\varphi \in\mathcal{D}(P(G))$  the
         integral $\int \d x \varphi(x) \langle f_1 | \kappa(x)$ exists and we have
         $$(\int \d x \varphi(x) \langle f_1 | \kappa(x))(f_2)=\int \d x \varphi(x) \langle f_1 | \kappa(x)|f_2 \rangle
         = \langle f_1| \mu(\varphi)| f_2 \rangle .$$
         Hence
         $$\int \d x \varphi(x) \langle f_1 | \kappa(x)=	\langle f_1| \mu(\varphi)\in L^2(G).$$
         Similarly define $\kappa(x)|f \rangle$ and obtain
         $$\int \d x \varphi(x) \kappa(x)|f \rangle=	 \mu(\varphi) f_2 \rangle \in L^2(G).$$
         So we can form the scalar product
         $$(\int \d x_1 \varphi_1(x_1) \langle f_1| \kappa(x_1)) (\int \d x_2 \varphi_2(x_2)  \kappa(x_2)| f_2 \rangle
         = \langle f_1| \mu(\varphi_1)\mu(\varphi_2)| f_2 \rangle
         =\langle f_1| \mu(\varphi_1\varphi_2)| f_2 \rangle $$
      This proves the first equation. We check the last equation
      \begin{multline*}
       \langle \alpha'(x)| \alpha(x')\rangle=  N(x)N(x')(C_1(x)\langle h| \delta(x-\Omega)+C_2(x)\langle h|\frac{\mathcal{P}}{x-\Omega})
       (C_1(x') \delta(x'-\Omega)|g \rangle+ C_2(x')\frac{\mathcal{P}}{x'-\Omega} |g \rangle)\\
       = N(x)N(x')\iint \d w_1 \d w_2 \langle h|\delta(w_1-\Omega)\delta(w_2-\Omega)| g \rangle \\
      (C_1(x)\delta(x-w_1) + C_2(x)\frac{\mathcal{P}}{x-w_1})
        (C_1(x')\delta(x'-w_2) + C_2(x)\frac{\mathcal{P}}{x'-w_2})\end{multline*}
        Use $\delta(w_1-\Omega)\delta(w_2-\Omega)=\delta(w_1-w_2)\delta(w_1-\Omega)$  and continue
        $$=N(x)N(x')\int \d w C_2(w)(C_1(x)\delta(x-w) + C_2(x)\frac{\mathcal{P}}{x-w})
        (C_1(x')\delta(x'-w) + C_2(x')\frac{\mathcal{P}}{x'-w}).$$
        Use the properties of the $ \delta$- function and equation (\ref{eq3})
    and continue
       \begin{multline*} = N(x)N(x')(C_2(x)C_1(x)^2 \delta(x-x')+C_2(x)C_2(x')(C_1(x')\frac{\mathcal{P}}{x-x'}+ C_1(x)\frac{\mathcal{P}}{x'-x})
       + \pi^2 C_2(x)^2 \delta(x-x'))\\
               + C_2(x)C_2(x')\int\ d w C_2(w)\frac{1}{x'-x}(\frac{\mathcal{P}}{x-w}- \frac{\mathcal{P}}{x'-w}))\end{multline*}
               Now
       \begin{multline*}
       C_1(x')\frac{\mathcal{P}}{x-x'}+ C_1(x)\frac{\mathcal{P}}{x'-x}=\frac{C_1(x')-C_1(x)}{x-x'}= \langle h| \frac{1}{x-x'}
       (-\frac{\mathcal{P}}{x'-\Omega}+\frac{\mathcal{P}}{x-\Omega})|g \rangle\\
       = \int \d w C_2(w)  \frac{1}{x-x'}
       (-\frac{\mathcal{P}}{x'-w}+\frac{\mathcal{P}}{x-w}))\end{multline*}
       The verification of the two other relations are left to the reader.
       \end{pf}
       \subsection{Perturbation of the Multiplication Operator 2}     This example is a caricature of the eigenvalue problem arising in the
        theory of radiation transfer
       in a gray atmosphere in plan parallel geometry \cite{EvW}.  We consider for some $c>1$ the set 
       $ G = ]-c,-1[ \cup ]1,c[\subset \mathbb{R}$, the multiplication operator
       $ \Omega f(y)=yf(y)$
       and two real $C^\infty $ functions $f,g$ on $\mathbb{R}$
        with $g(x)>0$ for $ 1<x<c$ and $-c < x< -1$ and 0 outside these two open interval. We assume  $ g(y)=g(-y)$ and
        $h(y)=-g(y)$ for $y>	1$ and $h(y)=g(y)$ for $y<-1$.   We study $H= \Omega + |g \rangle \langle h |$ and obtain
     $$ C(z)= 1 - \int_G\d y \frac{g(y)h(y)}{z- y}= 1 + \int _1^c \d y \frac{2 y g(y)^2}{ z^2-y^2}.$$
     We have for $z = x + \i u$
      $$ \Im  C(z) = \int _1^c \d y  \frac{4 y^2 xu g(y)^2}{(x^2-u^2 -y^2)^2+ 4 x^2u^2}$$
         Hence $C(z)=0$ implies $xu=0$ , so either $x$ or $u$ or both vanish. 
         We have
         \begin{align*}&                                                         
         C(0)=1-\int_1^c \d y \frac{2g(y)^2}{y},&&
       C(\i u)= 1 - \int_1^c \d y \frac{ 2 y g(y)^2}{u^2+y^2}\end{align*}
      So $C(\i u)$ is monotonic increasing for increasing $u^2$ and goes to 1 for $u^2 \to \infty$.
       If $ C(0) < 0$, there exists exactly one $ u_0> 0$
    such that $C(\i u_0)=0$, if $C(0) >0$, then $C(\i u )> 0$ for all $u$.
    
    For $| x| \le 1$ we have
    $$ C(x)=1-\int_1^c \d y \frac{2yg(y)^2}{y^2-x^2}$$
    and is monotonic decreasing for increasing $x$. If $C(0)>0$ and $C(1)<0$ there exists exactly one $x_0$ with $0<x_0<1$, such that
    $C(x_0)=0$. If $C(0)<0$ there does not exist such an $x$. We do not discuss the case $C(0)> 0$ and $C(1) \ge 0$. In case
    $C(0)=0$ we have a double zero. For $|x| \ge c$ we have $C(x) \ge 1$.
Hence $C(x)$ does not vanish for $|x| \ge c $.

The singularities of the resolvent  are the slits $[-c,-1]$ and $[1,c]$ and the zeros of $C(z)$. In the neighborhood of a zero   of $C(z)$
we may define a resolvent distribution with the help of proposition \ref{prop3}.
 We discuss  the behavior of the resolvent in the
neighborhood of 
the slits. We have
$$ C(x \pm \i 0)=C_1(x) \pm \i \pi C_2(x)= 1- \int \d y  g(y) h(y)\frac{\mathcal{P}}{x-y} \pm \i\pi  g(x)h(x). $$
There exists a neighborhood $G_1\subset \mathbb{R}$ of $[-1,-c]\cup [1,c]$ and an open interval $I$ containing $0$, such that $R(z)$ is
 holomorhic in $(G_1\times I) \setminus( [-1,-c]\cup [1,c])$ and $C(x \pm \i 0) \ne 0$ for $x \in G_1$. So we can apply prop. \ref{prop4}
  and define
 a distribution $R$ extending the resolvent function $R(z)$ to $G_1\times I$ and fulfilling the distribution resolvent equation
 . From the local definition in the neighborhood of the singularities we can define a resolvent distribution extending $R(z)$ to
 $\mathbb{C}$ with the help of proposition \ref{prop2}. 
       \begin{prop}
We calculate the spectral distribution. If there are two zeros $\ne 0$ we obtain
$$ M(z)=M(x+\i u)=r_+\delta(z-z_0)+r_-\delta(z+z_0)+\delta(u)\mu(x)$$
where $r_\pm$ are the residues of $R(z)$ at the points $\pm z_0$ and
$$ \mu(x)= \frac{1}{2\pi \i}(R(x-\i 0)-R(x+\i 0).$$
We identify $\mu(x)$ with its restriction $\kappa(x)$ as bilinear form on $\mathcal{D}(\mathbb{C})$ and obtain
$$ M(z) = |\alpha_+\rangle \langle\alpha_+'|\delta(z-z_0)+ |\alpha_- \rangle \langle \alpha_-'|\delta(z+z_0)+
 \delta(u) |\alpha_x \rangle \langle \alpha_x' |$$
 with the right resp. left usual eigenvectors
 \begin{align*}
&| \alpha_\pm \rangle = \frac{1}{\sqrt{ \langle h|(\pm z_0-\Omega)^{-2}|g \rangle}}\frac{1}{ \pm z_0-\Omega}| g \rangle&&
 \langle\alpha_\pm' |= \frac{1}{\sqrt{ \langle h|(\pm z_0-\Omega)^{-2}|g \rangle}}\langle h |\frac{1}{ \pm z_0-\Omega}
 \end{align*}
 and for $x \in G$ the right, resp. left generalized eigenvectors
 \begin{align*}
 | \alpha_x \rangle = \frac{1}{\sqrt{C_1^2 + \pi^2C_2^2}}(C_1(x)| \delta_x \rangle + h(x)A(x))&&
  \langle\alpha_x' |= \frac{1}{\sqrt{C_1^2 + \pi^2C_2^2}}(C_1(x) \langle\delta_x |  + g(x)A'(x))
  \end{align*}
  with
  \begin{align*}&
   C_1(x)==1-\int \d y g(y)h(y)\frac{\mathcal{P}}{x-y}&&
  C_2(x)=g(x)h(x)\\&
  A(x)=\frac{\mathcal{P}}{x-\Omega}| g \rangle &&
  A'(x)= \langle h |\frac{\mathcal{P}}{x-\Omega} \\&
  B(x) = g(x) | \delta_ x \rangle&&
  B'(x)=h(x) \langle \delta_x |\end{align*}
 In the case $C(0)=0$ we obtain
 $$R(z)=z^{-2}\frac{ \Omega^{-1}| g \rangle \langle h | \Omega^{-1}}{\langle h| \Omega^{-3}|g \rangle}+
 z^{-1} \frac{\Omega^{-2}| g \rangle \langle h | \Omega^{-1}+
\Omega^{-1}| g \rangle \langle h | \Omega^{-2}
}{\langle h| \Omega^{-3}|g \rangle}+ O(1)= z^{-2}a+ z^{-1}p_0+ O(1)$$
and
$$M(z)=M(x + \i u)= \ p_0 \delta_2(z) - a \partial  \delta_2(z) + \mu(x)\delta(u),$$  
where $\mu(x) =| \alpha_x 
\rangle \langle\alpha'_x | $ is given by the formula above. 
Here $p^2=p$ and $a^2=0$ and $ap=pa=a$. 
We obtain for $\vartheta,\vartheta' = \pm$ and $x,x' \in \mathbb{R}$ the orthogonality relations
\begin{align*}&
 \langle \alpha_\vartheta' | \alpha_{\vartheta'}\rangle = \delta_{\vartheta,\vartheta'}&&
 \langle \alpha_x'| \alpha_\vartheta \rangle =\langle \alpha_\vartheta | \alpha_ x\rangle = 0 &&
\langle \alpha_x| \alpha_{x'}\rangle = \delta(x-x')\end{align*}
Analogous relations hold for the case $C(0)=0$. The spectral distribution is complete, i.e.
$$M(1)= \int ]d^2 z M(z)= \sum_{\vartheta = \pm}| \alpha_\vartheta ' \rangle \langle \alpha _\vartheta |
+ \int_g \d x | \alpha_x ' \rangle \langle \alpha _x |= 1$$
\end{prop}
\begin{pf}
Assume at first, that $C(z)$ has two zeros $\pm z_0= \pm \i u_0$ or $ \pm z_0= \pm x_0$.
 The residuum of
the complex function $R(z)$ at the points $ \pm z_0$ is given by
$$r_\pm= \big(\frac{1}{\pm z_0 - \Omega}|g \rangle \langle h | \frac{1}{\pm z_0-\Omega}\big)/
C'(\pm z_0)= \big(\frac{1}{\pm z_0 - \Omega}|g \rangle \langle h | \frac{1}{\pm z_0-\Omega}\big)
\frac{1}
{\langle h| (\pm z_0-\Omega)^{-2}|g \rangle}$$
Using the results of  proposition \ref{prop3} we obtain
$$ M(z) = M(x+\i u)= r_+\delta_2(z-z_0)+ r_-\delta_2(z+z_0)+ \delta(u)\mu(x)$$
where 
      $$\mu(x)= \frac{1}{2 \pi \i}(R(x-\i 0)-R(x+ \i 0)$$
   We  consider the case $C(0)=0$ .We expand at the origin
$$ C(z)= 1+ \sum_{n=0}^\infty z^n \langle h| \Omega ^{-(n+1)}|g \rangle =
z^2 \langle h| \Omega^{-3}|g \rangle+ O(z^3)=
-z^2 \int _1^\infty 2g(y)/y^3 \d y + O(z^3)  $$

$$R(z)=z^{-2}\frac{ \Omega^{-1}| g \rangle \langle h | \Omega^{-1}}{\langle h| \Omega^{-3}|g \rangle}+
 z^{-1} \frac{\Omega^{-2}| g \rangle \langle h | \Omega^{-1}+
\Omega^{-1}| g \rangle \langle h | \Omega^{-2}
}{\langle h| \Omega^{-3}|g \rangle}+ O(1)= z^{-2}a+ z^{-1}p_0+ O(1)$$
and
$$M(z)=M(x + \i u)= \ p_0 \delta_2(z) - a \partial  \delta_2(z) + \mu(x)\delta(u).$$  
Here $p^2=p$ and $a^2=0$ and $ap=pa=a$. 

In both cases we can use for $f_1,f_2 \in \mathcal{D}(\mathbb{R})$ the relation  $\langle f_1| \mu(x)| f_2 \rangle =
\langle f_1| \kappa(x)| f_2 \rangle$ and $\kappa(x)$ is given by the formula   
\begin{multline*} 
         \kappa(x) = \delta(x-\Omega)+ \frac{1}{C_1^2+\pi^2 C_2^2}
       (AC_2A'+AC_1B'+BC_1A'- \pi^2 BC_2B')   \\ 
 = \frac{1}{C_1^2+\pi^2 C_2^2}(Ah(x)g(x)A'+AC_1h(x) \langle \delta_x|
+ g(x)|\delta_x \rangle C_1 A' + C_1^2 h(x)^2 g(x)^2 | \delta_x \rangle \langle \delta_x |)= | \alpha_x 
\rangle \langle\alpha'_x | \end{multline*}
as $\delta(x-\Omega)=|\delta_x \rangle \langle \delta_x|$.
The proof of the orthogonality relations can be deduced from  the relation $\kappa(x)\kappa(x')= \delta(x-x')\kappa(x)$, like in the proof of 
prop  \ref{prop8} or can be done by hand. That
 is  left to the reader. For  the completeness refer to theorem \ref{thm2}.
 \end{pf}

\end{document}